\documentclass[11pt,a4paper]{amsart}

\usepackage{amsmath}
\usepackage{amssymb}
\usepackage{caption2}
\usepackage[tight,TABTOPCAP]{subfigure}
\usepackage{tabularx}
\usepackage{hyperref}

\makeatletter
\renewcommand{\@thesubtable}{\hspace*{-4em}}
\makeatother

\numberwithin{equation}{section}

\newcommand{\F}{\mathbb{F}_2}
\renewcommand{\>}{\rangle}
\newcommand{\<}{\langle}
\newcommand{\su}{{\rm supp}}
\newcommand{\End}{\operatorname{End}}
\newcommand{\card}[1]{{\mid\! #1 \!\mid}}
\newcommand{\pro}[2]{\langle #1, #2 \rangle}

\def\reg{{\rm reg}}
\def\lra{{\,\Leftrightarrow\,}}
\def\ra{{\,\Rightarrow\,}}
\def\Irr{{\rm Irr}}
\def\C{{\mathbb C}}
\def\N{{\mathbb N}}
\def\R{{\mathbb R}}
\def\Z{{\mathbb Z}}

\def\V{{\mathcal V}}
\def\conv{{\rm conv}}
\def\aff{{\rm aff \,}}
\def\lin{{\rm lin \,}}
\def\dim{{\rm dim \,}}
\def\deg{{\rm deg \,}}
\def\lcm{{\rm lcm}}
\def\GL{{\rm GL}}
\def\Mat{\textup{Mat}}
\def\id{\textup{id}}

\newtheorem*{theorem*}{Theorem}
\newtheorem*{corollary*}{Corollary}
\newtheorem{theorem}{Theorem}[section]
\newtheorem{definition}[theorem]{Definition}
\newtheorem{remark}[theorem]{Remark}
\newtheorem{example}[theorem]{Example}
\newtheorem{question}[theorem]{Question}
\newtheorem{corollary}[theorem]{Corollary}
\newtheorem{prop}[theorem]{Proposition}
\newtheorem{lemma}[theorem]{Lemma}
\newtheorem{conjecture}[theorem]{Conjecture}
\newtheorem*{acknowledgement}{Acknowledgement}

\begin{document}
\title{On permutation polytopes}

\author[Baumeister]{Barbara Baumeister}
\email{\{baumeist,christian.haase,nill,paffenho\}@math.fu-berlin.de}
\author[Haase]{Christian Haase}
\author[Nill]{Benjamin Nill}
\author[Paffenholz]{Andreas Paffenholz}
\address{Mathematisches Institut, Freie Universit\"at Berlin,
  Arnimallee 3, 14195 Berlin, Germany}

\begin{abstract}
  A permutation polytope is the convex hull of a group of
  permutation matrices. In this paper we investigate the combinatorics
  of permutation polytopes and their faces. As applications we
  completely classify $\leq 4$-dimensional permutation polytopes and
  the corresponding permutation groups up to a suitable notion of equivalence. 
  We also provide a list of combinatorial types 
  of possibly occuring $\leq 4$-faces of permutation polytopes.
\vspace{-2ex}
\end{abstract}

\maketitle

\section*{Introduction}
One of the most intensively studied convex polytopes is the Birkhoff
polytope, also known as the assignment polytope, also known as the
polytope of doubly stochastic matrices \cite{BR74,BG77,BL91,BS96,Zei99,CRY00,Pak00}. 
It is the convex hull in $\R^{n \times n}$ of the $n \times n$
permutation matrices. This polytope naturally appears in various contexts 
such as enumerative combinatorics \cite{Sta86,Ath05},
optimization \cite{Tin86,Fie88,Pak00,BS03}, and statistics \cite{Pak00} 
(and references therein), as well as in 
representation theory \cite{Onn93,BFL+02}, and 
in the context of the van der Waerden conjecture for the
permanent \cite{BG77}.

In the present article, we propose to systematically study 
general permutation polytopes. These are defined as the convex hull of a subgroup $G$
of the group of $n \times n$ permutation matrices. This is a convex geometric
invariant of a permutation representation, and it yields various
numerical invariants like dimension, volume, diameter, $f$-vector,
etc.

A number of authors have studied special classes of permutation
polytopes different from the Birkhoff polytope. Brualdi and Liu~\cite{BL91} compute basic invariants of the
polytope of the alternating group; for this polytope, Hood and Perkinson~\cite{JH04}
describe exponentially many facets.
Collins and Perkinson~\cite{CP04} observe that Frobenius polytopes
have a particularly simple combinatorial structure, and 
Steinkamp~\cite{Ste99} adds results about dihedral groups. Most recently, 
Guralnick and Perkinson~\cite{GP06} investigate general permutation
polytopes, their dimension, and their graph.

\subsection*{Main results}
In Section~\ref{sec:notation1}, we introduce the main objects of our study, 
representation polytopes and permutation polytopes. We also add a 
note on $0/1$-polytopes, Proposition \ref{01bound}, showing that there is only a finite number of 
lattice equivalence classes of $d$-dimensional $0/1$-polytopes.

In Section~\ref{sec:notation2}, we discuss notions of equivalence of 
representations and the associated representation polytopes, respectively permutation polytopes. 
In particular, we introduce stable and effective equivalence of representations.

In Section~\ref{sec:combinatorics}, we investigate combinatorial
properties of permutation polytopes and their faces. In particular, we are interested in the question 
which polytopes can be realized as faces of permutation polytopes. The first main theorem, Theorem \ref{product}, 
says that, if a permutation polytope is combinatorially a product, 
then the permutation group has a natural 
product structure. As a second result, we give an explicit construction, Theorem \ref{pyramid}, 
showing that pyramids over faces of permutation polytopes appear again as faces of permutation polytopes. 
Further, we are interested in centrally symmetric faces and polytopes. We show in Theorem \ref{constr} that 
free sums of crosspolytopes and cubes occur as faces of permutation polytopes. Finally, 
we construct in Theorem \ref{cross} the essentially unique permutation group where the permutation polytope 
is a crosspolytope.

In Section~\ref{sec:classification}, we use the results from the previous section to classify in Theorem \ref{permpolytheo} 
up to effective equivalence all permutation representations whose polytopes have dimension $\le 4$. 
We also start the more difficult classification of combinatorial
types of low dimensional polytopes which appear as faces of permutation polytopes. 
Theorem \ref{facestheo} gives a complete answer for dimension $\le 3$, and 
there remain only two $4$-polytopes for which we could not decide whether or not they can be realized. 
The lists of permutation groups and polytopes can also be found on the webpage \cite{BHNP07}.

Those examples and classifications suggest a number of open questions
and conjectures which we formulate in Section~\ref{sec:questions}. 

\subsubsection*{Remark:} Some authors use the notion of permutation polytope differently: the convex
hull of a $G$-orbit in $\R^n$. These are linear projections of our permutation
polytopes. Examples include the permutahedron, the traveling salesman
polytope, or any polytope with a vertex transitive group of
automorphisms. They appear in combinatorial optimization problems
of various computational complexities~\cite{Onn93}. Moreover, orbit
polytopes have recently been used to construct resolutions in group
cohomology~\cite{EHS06}.

\section{Representation-, permutation-, and $0/1$-polytopes}
\label{sec:notation1}
\subsection{Polytopes}

For a standard reference on polytopes we refer to \cite{Zie95}. 
A polytope $P$ is the convex hull $\conv(S)$ of a finite set of points
$S$ in a real vector space $V$. The dimension $\dim P$ is the
dimension of the affine hull $\aff P$ as an affine space. 
We say $P$ is a $d$-polytope, if $\dim P = d$. 
If $V$ is equipped with a full dimensional lattice $\Lambda$ and we
can choose $S \subset \Lambda$, then we call $P$ a lattice polytope.

A face $F$ of $P$ (denoted by $F \preceq P$) 
is a subset where some linear functional is maximized.
Zero-dimensional faces are vertices, one-dimensional faces are edges,
and faces of codimension one are facets. 
The poset of faces ordered by inclusion is called the face lattice. 
The vertex set of $P$ is denoted by $\V(P)$. The degree of a vertex is the number 
of edges it is contained in.

There is a hierarchy of equivalence relations on (lattice) polytopes. 
Two polytopes $P \subset \R^m$ and $Q \subset \R^n$ are affinely
equivalent if there is an affine isomorphism of the affine hulls $\phi
\colon \aff P \to \aff Q$ that maps $P$ onto $Q$. For lattice
equivalence we additionally require that $\phi$ is an isomorphism of
the affine lattices $\aff P \cap \Lambda \to \aff Q \cap
\Lambda'$. Combinatorial equivalence is merely an equivalence of the face
lattices as posets. 
\newlength{\widthOfEquivalent}
\settowidth{\widthOfEquivalent}{equivalent}
\newlength{\widthOfCombinatorially}
\settowidth{\widthOfCombinatorially}{combinatorially}
$$\raisebox{.1\baselineskip}{\parbox{\widthOfEquivalent}{lattice
    \\[-1mm] equivalent}}
\quad \Rightarrow \quad
\raisebox{.1\baselineskip}{\parbox{\widthOfEquivalent}{affinely
    \\[-1mm] equivalent}}
\quad \Rightarrow \quad 
\raisebox{.1\baselineskip}{\parbox{\widthOfCombinatorially}{combinatorially
    \\[-1mm] equivalent}}$$
The converse implications do not hold, for examples see \cite[Prop.~7]{Zie00}.

\subsection{Representation polytopes}

Let $\rho \colon G \to \GL(V)$ be a real representation of the finite
group $G$ with identity element $e$. It induces an $\R$-algebra homomorphism from the group algebra $\R[G]$ to $\End(V)$, 
which we also denote by $\rho$.

\begin{definition}{\rm 
The {\em representation polytope} $P(\rho)$ of the representation $\rho$ 
is defined as the convex hull of $\rho(G)$ in the vector space $\End(V)$.
}
\end{definition}

Notice, that the representation $\rho$ splits as a $G$-representation over $\C$ into irreducible components: 
\begin{equation}
\rho \cong \sum_{\sigma \in \Irr(G)} c_\sigma \sigma
\label{irrcomp}
\end{equation}
with $c_\sigma \in \N$ for $\sigma$ in $\Irr(G)$, the set of pairwise non-isomorphic 
irreducible $\C$-representations. We define the set of {\em irreducible factors} of $\rho$, 
\[\Irr(\rho) := \{\sigma \in \Irr(G) \,:\, c_\sigma > 0\}.\]

The group $G$ acts on the polytope $P(\rho)$ by left multiplication, inducing an affine automorphism of $P(\rho)$: 
$$g(\rho(h)) = \rho(g)\rho(h) \mbox{ for all }g,h \in G.$$

Therefore, since any vertex of $P(\rho)$ has to be contained in $\rho(G)$, and left multiplication on $G$ is regular, thus transitive, 
we get:
\begin{equation}\label{ecken}\V(P(\rho)) = \rho(G).\end{equation}
Here is one application (the case of equality is treated in Corollary \ref{simplex}):
\begin{equation}\label{upperdimension}
\dim P(\rho) \leq \card{\V(P(\rho))} - 1 \leq \card{G} - 1.\end{equation}
In particular, though there are infinitely many representations of finite groups of fixed order, 
they give rise to only finitely many combinatorial types of representation polytopes. 
We are going to see a stronger statement in Corollary \ref{finiteness}.

More implications: All vertices of 
$P(\rho)$ have the same degree. 
When considering the combinatorics of
a face of $P(\rho)$ we can always assume that it has $\rho(e)=\id$ as a vertex.
If $F$ is a face of $P(\rho)$ with vertex set $\V(F) = \rho(U)$ for $U \subseteq G$, then 
$F$ is also a face of the representation polytope $P(\rho')$, where $\rho' \colon \<U\> \to \GL(V)$. 
Here $\<U\>$ denotes the smallest subgroup of $G$ containing~$U$.

\subsection{Permutation polytopes}

We identify the symmetric group $S_n$ on $\{1, \ldots, n\}$ via the usual
permutation representation with the set of $n \times n$ permutation
matrices, i.e., the set of matrices with entries $0$ or $1$ such that
in any column and any row there is precisely one $1$. Throughout, we use cycle notation: 
For instance $(1 2 3) (4 5) \in S_6$ denotes the permutation $1 \mapsto 2$, $2 \mapsto~3$, $3 \mapsto 1$, $4 \mapsto 5$, 
$5 \mapsto 4$, $6 \mapsto 6$. Note that for $g_1, g_2 \in S_n$ we have $g_1 g_2 \in S_n$, while
$g_1 + g_2 \in \Mat_n(\N)$. Here, for a set $C \subseteq \C$, we define $\Mat_n(C)$ as the set of $n
\times n$ matrices with entries in $C$. We identify $\Mat_n(\R) \cong \R^{n^2}$, thus 
we have the usual scalar product, i.e., $\pro{A}{B} = \sum_{i,j} A_{i,j} B_{i,j}$ for $A,B \in \Mat_n(\R)$.

A subgroup of $S_n$ is called permutation group. 
A faithful representation $\rho \colon G \to S_n$ is called permutation representation, thus 
$G$ can be identified with the permutation group $\rho(G)$. For both situations 
we often write in short $G \leq S_n$.

\begin{definition}{\rm For $G \leq S_n$ we define $P(G) := \conv(G) \subseteq \Mat_n(\R)$, the {\em permutation polytope} associated to $G$. 
The convex hull of all permutations $B_n := P(S_n)$ is called the $n^{th}$ {\em Birkhoff polytope}.
}\end{definition}

In particular, any permutation polytope is a representation polytope, as well as a 
lattice polytope with respect to the lattice $\Mat_n(\Z)$.

\subsection{$\mathbf{0/1}$-polytopes}

An important property of permutation polytopes is that they belong
to the class of $0/1$-polytopes. A $0/1$-polytope is the convex hull of
points in $\{0,1\}^d$. 
They have been classified up to dimension $6$ \cite{Aic00,Aic07}. For a survey on these well-studied polytopes see
\cite{Zie00}, where also the following basic fact is shown:  any
$d$-dimensional $0/1$-polytope is affinely equivalent to a
lattice polytope in $[0,1]^d$.
This implies immediately that there are only finitely many affine
types of $d$-dimensional $0/1$-polytopes. Even more is true.

\begin{prop}
  Every $d$-dimensional $0/1$-polytope is lattice equivalent to a
  lattice polytope in the $2^{2^d}$-dimensional unit cube. In particular, there are only
  finitely many lattice types of $d$-dimensional $0/1$-polytopes.
\label{01bound}
\end{prop}
This bound is far from optimal. All we care about is that it is
finite. We have not found this result in the literature, so we
include the proof.
\begin{proof}
  Every $d$-dimensional $0/1$-polytope is lattice equivalent to a
  full-di\-men\-sional lattice polytope which will, in general, no
  longer have $0/1$ coordinates. But it still has the property that
  its vertices are the only lattice points it contains. Such a
  polytope can have no more than $2^d$ vertices, as two vertices with
  the same parity would have an integral midpoint.

  Now suppose $P \subset [0,1]^N$ is a $d$-dimensional $0/1$-polytope
  with $N > 2^{2^d}$. Then there are two of the $N$ coordinates which
  agree for every vertex of $P$, say, $P \subset \{ x_i=x_j \}$. Thus
  we can delete the $j^{\mathrm th}$ coordinate and obtain a lattice
  equivalent $0/1$-polytope $\subset [0,1]^{N-1}$.
\end{proof}

\begin{corollary}
There are up to lattice equivalence only finitely many permutation polytopes associated to finite groups of fixed order.
\label{finiteness}
\end{corollary}

This follows from Equation (\ref{upperdimension}).

\section{Notions of equivalence}
\label{sec:notation2}

Throughout let $\rho \colon G \to \GL(V)$ be a real representation. 

\subsection{Stable equivalence of representations}

When working with permutation polytopes, one would like to identify representations that 
define affinely equivalent polytopes. For instance, this holds for the following five permutation groups: 
$\<(1 2 3 4)\> \leq S_4$, $\<(1 2 3 4) (5)\> \leq S_5$, 
$\<(1 2 3 4) (5 6)\> \leq S_6$, $\<(1 2 3 4) (5 6) (7 8)\> \leq S_8$, $\<(1 2 3 4) (5 6 7 8)\> \leq S_8$. 
We are now going to introduce a suitable notion of equivalence on the real representations of a finite group. 
The crucial point is the observation that representation polytopes do not care about multiplicities 
of irreducible factors in the defining representation. For this let us fix a finite group $G$.

\begin{definition}{\rm 
  For a representation $\rho \colon G \to \GL(V)$ define the affine
  kernel $\ker^\circ \rho$ as
  $$\ker^\circ \rho := \left\{ \ \sum_{g \in G} \lambda_g g \in \R[G] \ : \
    \sum_{g \in G} \lambda_g \rho(g) = 0 \text{ and } \sum_{g \in G} \lambda_g = 0
    \ \right\}$$

  Say that a real representation $\rho' \colon G \to \GL(V')$ is an {\em affine
  quotient} of $\rho$ if $\ker^\circ \rho \subseteq \ker^\circ \rho'$ .

  Then real representations $\rho_1$ and $\rho_2$ of $G$ are {\em stably equivalent}, if
  there are affine quotients $\rho_1'$ of $\rho_1$ and $\rho_2'$ of
  $\rho_2$ such that $\rho_1 \oplus \rho_1' \cong \rho_2 \oplus
  \rho_2'$ as $G$-representations. For instance, $\rho_1$ is stably equivalent to $\rho_1 \oplus \rho'_1$.
}
\end{definition}

\begin{example}{\rm \label{ex:affine quotient}
Let $1_G$ be the trivial representation of $G$. We observe that $\ker^\circ \rho = \ker 1_G \cap \ker \rho$. Hence, by 
Equation (\ref{irrcomp})
\[\ker^\circ \rho = \ker 1_G \;\;\cap \bigcap_{1_G \not= \sigma \in \Irr(\rho)} \ker \sigma.\]
Therefore, 
any real representation $\rho' \colon G \to \GL(V')$ with $\Irr(\rho')\backslash\{1_G\} \subseteq \Irr(\rho)\backslash\{1_G\}$ 
is an affine quotient of $\rho$. For instance, $\rho'$ may be the restriction of $\rho$ to an invariant subspace of $V$.
}
\end{example}

\begin{prop} \label{prop:stablyEquivalence}
  Suppose $\rho$ and $\bar{\rho}$ are stably equivalent real representations of a finite group $G$. 
Then $P(\rho)$ and $P(\bar{\rho})$ are affinely equivalent.
\end{prop}

\begin{proof}
  It is enough to show 
that $P(\rho)$ and $P(\rho \oplus \rho')$ are affinely equivalent for an affine quotient $\rho'$ of $\rho$.

  The projection yields an affine map $P(\rho \oplus \rho') \to
  P(\rho)$. In order to construct an inverse, we need a map $\aff P(\rho) \to
  \aff P(\rho')$. The obvious choice is to map a point $\sum_{g \in G}
  \lambda_g \rho(g) \in \aff P(\rho)$ to $\sum_{g \in G} \lambda_g \rho'(g) \in
  \aff P(\rho')$ ($\sum_{g \in G} \lambda_g = 1$). This is well defined if
  (and only if) $\ker^\circ \rho \subseteq \ker^\circ \rho'$. 
\end{proof}

A priori, it is often not clear whether two representations are stably 
equivalent. Here we provide an explicit criterion:

\begin{theorem}\label{thm:stablyEquivalence}
Two real representations are stably equivalent if and only if they contain 
the same non-trivial irreducible factors.
\end{theorem}

The proof will be given in the next subsection. 
We note that we have already seen the if-direction in Example \ref{ex:affine quotient}.

\subsection{The dimension formula}

To prove Theorem \ref{thm:stablyEquivalence} we recall the dimension formula of a 
representation polytope in \cite{GP06}.

The following equation is Theorem 3.2 of \cite{GP06} (recall that the degree of a representation is the dimension of 
the vector space the group is acting on).

\begin{theorem}[Guralnick, Perkinson]
\[\dim P(\rho) = \sum_{1_G \not= \sigma \in \Irr(\rho)} (\deg \sigma)^2.\]
\label{dimension}
\end{theorem}

The proof relies on the Theorem of Frobenius and Schur \cite[(27.8-10)]{CR62} to determine the dimension of $\rho(\C[G])$, 
and then relates $\dim P(\rho)$ to $\dim_\C \rho(\C[G])$ via the following observation which explains the special 
role of the trivial representation.

\begin{lemma}\label{lem:zero vector}
The affine hull of $P(\rho)$ does not contain $0$ if and only if $1_G \in \Irr(\rho)$.
\end{lemma}

Now, we can give the proof of the characterization of stable equivalence:

\begin{proof}[Proof of Theorem \ref{thm:stablyEquivalence}]

It is enough to show 
that $\rho$ and $\rho \oplus \rho'$ have the same non-trivial irreducible factors for an affine quotient $\rho'$ of $\rho$.

By Proposition \ref{prop:stablyEquivalence} $P(\rho)$ and $P(\rho\oplus\rho')$ are affinely equivalent. 
In particular, they have the same dimension. Since any irreducible factor of $\rho$ is an irreducible factor of $\rho \oplus \rho'$, 
the dimension formula Theorem \ref{dimension} implies that any non-trivial irreducible factor of $\rho'$ already appears as 
an irreducible factor of $\rho$. Therefore, $\rho$ and $\rho'$ have the same non-trivial irreducible factors.
\end{proof}

For an application let us look at the regular representation of a group $G$. 
This is the permutation representation $\reg \colon G \to S_{\card{G}}$ 
via right multiplication. We have $\Irr(\reg) = \Irr(G)$.

\begin{lemma}
$P(\reg)$ is a simplex of dimension $\card{G} - 1$, and the vertices 
form a lattice basis of the lattice $\lin P(\reg) \cap \Mat_{\card{G}}(\Z)$.
\label{regular}
\end{lemma}

\begin{proof}

For this we enumerate the elements of $G$ as $g_1, \ldots, g_{\card{G}}$ with $g_1 = e$. 
Then for $i \in\{1, \ldots, \card{G}\}$ the permutation matrix $\rho(g_i)$ of size $\card{G} \times \card{G}$ has in the first row 
only zeros except one $1$ in column $i$. Hence, the matrices $g_1, \ldots, g_{\card{G}}$ are linearly independent. Moreover, this shows that 
they form a lattice basis of $\lin P(\reg) \cap \Mat_{\card{G}}(\Z)$. 
\end{proof}

The dimension formula and Lemma \ref{regular} imply another proof of the following well-known equation:
\[\card{G}-1 = \dim P(\reg) = \sum_{1_G \not= \sigma \in \Irr(G)} (\deg \sigma)^2.\]

We see that in a special case there is indeed a correspondence between stable equivalence and affine equivalence, 
this is \cite[Cor.~3.3]{GP06}.

\begin{corollary}
Let $\rho$ be a faithful representation. Then 
$P(\rho)$ is a simplex if and only if $\rho$ is stably equivalent to $\reg$.
\label{simplex}
\end{corollary}

Here is an example showing that stably equivalent permutation representations do not necessarily 
have lattice equivalent permutation polytopes.

\begin{example}{\rm 
Let $G:= \<(1 2),(3 4)\> \leq S_4$. We define the following permutation representation: 
$\rho \colon G \to S_6$, by $(1 2) \mapsto (1 2)(3 4)$ and $(3 4) \mapsto (1 2)(5 6)$. Then $P(\rho)$ is a tetrahedron, and $\rho$ is 
stably equivalent to the regular representation. However, the vertices of $P(\rho)$ do {\em not} form an affine lattice basis of the lattice 
$\aff P(\reg) \cap \Mat_{\card{G}}(\Z)$, in contrast to $P(\reg)$ by Lemma \ref{regular}.
}
\end{example}

Affine equivalence is the same as lattice equivalence for the sublattice generated by the vertices. 
The previous example shows that lattice equivalence for the whole lattice $\Mat_n(\Z)$ is a more subtle condition. 
This relation deserves further study.

\subsection{Effective equivalence of representations}

The following example illustrates that stable equivalence is too rigid.

\begin{example}{\rm 
Let $G:= \<(1 2),(3 4)\> \leq S_4$. Then $\rho_1 \colon G \hookrightarrow S_4$ is a permutation representation with 
$P(G) = P(\rho_1)$ a square. On the other hand, we define another permutation representation 
$\rho_2 \colon G \to S_4$, by $(1 2) \mapsto (1 2)$ and $(3 4) \mapsto (1 2)(3 4)$. Then $P(\rho_2)$ is the same 
square. However, $\rho_1$ and $\rho_2$ are {\em not} stably equivalent, since they do not have 
the same irreducible factors. 
\label{vier}
}
\end{example}

We observe that these two representations $\rho_1, \rho_2 \colon G \to \GL(V)$ are conjugated, i.e., 
there exists an automorphism $\psi$ of $G$ such that $\rho_2 = \rho_1 \circ \psi$. Hence, since $\rho_1(G) = \rho_2(G)$, 
we have $P(\rho_1) = P(\rho_2)$. However, conjugation permutes the irreducible factors, thus does not respect stable equivalence. 
To avoid this ambiguity we propose the following notion.

\begin{definition}{\rm Two real representations $\rho_i \colon G_i \to \GL(V_i)$ (for $i=1,2$) of finite groups 
are {\em effectively equivalent}, if there exists an isomorphism $\phi \colon G_1 \to G_2$ such that 
$\rho_1$ and $\rho_2 \circ \phi$ are stably equivalent $G_1$-representations. 

Moreover, we say $G_1 \leq S_{n_1}$ and $G_2 \leq S_{n_2}$ are {\em effectively equivalent} permutation groups, if 
$G_1 \hookrightarrow S_{n_1}$ and $G_2 \hookrightarrow S_{n_2}$ are effectively equivalent permutation representations.
}
\end{definition}

By Theorem \ref{thm:stablyEquivalence} we may put this definition in a nutshell: Two permutation groups are effectively equivalent 
if they are isomorphic as abstract groups such that via this isomorphism the permutation representations 
contain the same non-trivial irreducible factors.

In particular, effectively equivalent representations have affinely equivalent representation polytopes by Proposition \ref{prop:stablyEquivalence}. 
Of course, in general the converse cannot hold, since by Lemma \ref{regular} permutation groups that are not even isomorphic as abstract groups 
still may have affinely equivalent permutation polytopes. However, the following question remains open.

\begin{question}{\rm Are there permutation groups $G_1, G_2$ that are isomorphic as abstract groups and whose 
permutation polytopes $P(G_1)$ and $ P(G_2)$ are affinely equivalent, while $G_1$ and $G_2$ are not effectively equivalent?
}
\end{question}

By Theorem \ref{permpolytheo} there are no such permutation groups, if their permutation polytopes have dimension $\leq 4$.

\section{The combinatorics of permutation polytopes} \label{sec:combinatorics}

Throughout, $G \leq S_n$ is a permutation group. By $G \cong H$ we denote an (abstract) group isomorphism.

\subsection{The smallest face containing a pair of vertices}\

In \cite{BG77,BL91} the diameter of the edge-graph of $B_n$ and $P(A_n)$ was bounded from above by $2$. 
Later this could be generalized in \cite{GP06} 
to  permutation polytopes associated to transitive permutation groups. 
For this, Guralnick and Perkinson needed a crucial observation, that we are going to recall here. 

\begin{definition}{\rm 
Let $e \not= g \in G$.
\begin{itemize}
\item The support $\su(g)$ is the complement of the set of fixed points.
\item We denote by $F_g$ the smallest face of $P(G)$ containing $e$ and $g$.
\item We denote by $g = z_1 \circ \cdots \circ z_r$ the unique disjoint cycle decomposition of 
$g$ in $S_n$, i.e., $z_1, \ldots, z_r$ are cycles with pairwise 
disjoint support, and $g = z_1 \cdots z_r$.
\item Let $g = z_1 \circ \cdots \circ z_r$. For $h \in S_n$ we say $h$ is a {\em subelement} of $g$, 
if there is a subset $I \subseteq \{1, \ldots, r\}$ such that $h = \prod_{i \in I} z_i$.
\item $g$ is called {\em indecomposable} in $G$, if $e$ and $g$ are the only subelements of $g$ in $G$.
\end{itemize}}
\label{indecomp}
\end{definition}

The following result is Theorem 3.5 in \cite{GP06}. We include the very instructive proof here.

\begin{theorem}[Guralnick, Perkinson]
Let $g \in G$. The vertices of $F_g$ are precisely the subelements of $g$ in $G$. In particular, 
$e$ and $g$ form an edge of $P(G)$ if and only if $g$ is indecomposable in $G$.
\label{smallfacetheo}
\end{theorem}

Therefore, the number of indecomposable elements (different from $e$) in $G$ equals the degree of any vertex of $P(G)$.

\begin{remark}{\rm 
The proof of Guralnick and Perkinson uses a simple but effective way of defining certain faces of a permutation polytope $P(G)$ for $G \leq S_n$. 
These faces are the intersections of $P(G)$ with faces of the Birkhoff polytope $P(S_n)$. 
Since this method will also be used for several results of this paper, we give here the explicit description. 

Let $S \subseteq S_n$ be a subset of permutation matrices. 
We define the $n \times n$-matrix $M(S) := \max(\sigma \,:\, \sigma \in S)$, where the maximum is applied for any entry. Then 
$M(S)$ has only entries in $\{0,1\}$, thus 
$\pro{M(S)}{g} \leq n$ for any $g \in G$. Therefore, $F(S) := \{g \in G \,:\, \pro{M(S)}{g} = n\}$ is a face of $P(G)$. 
If $S \subseteq G$, then $S \subseteq \V(F(S))$.

{\em If $S \subseteq G$ and $\card{S} \leq 2$, then $F(S)$ is even the smallest face of $P(G)$ containing $S$.} 
This is part of the proof of Theorem \ref{smallfacetheo}.

While in the case of Birkhoff polytopes this implication holds also for $\card{S} \geq 3$, it is important to note that 
in the case of general permutation polytopes it usually fails. The following example illustrates this 
phenomenon. 
Let $z_1 := (1 2)$, $z_2 := (3 4)$, $z_3 := (5 6)$, $z_4 := (7 8)$. We define $G := \<z_1 z_2, z_1 z_3, z_1 z_4\> \leq S_8$, and 
$S := \{e, z_1 z_2, z_1 z_3\}$. Then $P(G)$ is a four-dimensional crosspolytope, i.e., the dual is a $4$-cube, 
and $P(G)$ contains a face with vertices $S$. However, $F(S)$ also contains the vertex $z_2 z_3$, so it is not the smallest face of $P(G)$ 
containing $S$.
\label{reduction}
}
\end{remark}

\begin{proof}[Proof of Theorem \ref{smallfacetheo}]

Let $S := \{e,g\} \subseteq G$, 
and $F(S)$ the face of $P(G)$ as defined in the previous remark. Then the vertices of $F(S)$ are precisely 
the subelements of $g$ in $G$. On the other hand, let $g = z_1 \circ \cdots \circ z_r$, and 
$h = z_1 \circ \cdots \circ z_s$ ($s \leq r$) be a subelement of $g$ in $G$. 
Then $h' := g h^{-1} = z_{s+1} \circ \cdots \circ z_r$ is also a subelement of $g$ in $G$. Now, the following identity of matrices holds:
\[e+g=h+h'.\]
Therefore, $F(S)$ is centrally symmetric with center $(e+g)/2$. Hence, $F(S)$ is 
the smallest face $F_g$ of $P(G)$ containing $S = \{e,g\}$.
\end{proof}

In particular we see from the proof that, if $g \in G$ and $h \in \V(F_g)$, 
then the antipodal vertex of $h$ in the centrally symmetric face $F_g$ is given by $g h^{-1}$ 
with $\su(h) \cap \su(g h^{-1}) = \emptyset$. 
Let us note this important restriction on the combinatorics of a permutation polytope.

\begin{corollary}
The smallest face containing a given pair of vertices of a permutation polytope is centrally symmetric.
\label{smallfacecoro}
\end{corollary}

This generalizes the well-known fact (e.g., see \cite[Thm.~2.5]{BS96}) that the smallest face of the Birkhoff polytope 
containing a pair of vertices is a cube. This strong statement is not true for general permutation polytopes. 
For instance in Corollary \ref{crosscoro} we show that crosspolytopes appear as faces of permutation polytopes.

\subsection{Products}

Products of permutation polytopes are again permutation polytopes, and therefore also 
products of faces of permutation polytopes appear as faces of permutation polytopes.

In many cases, given a permutation group $G$ and its permutation polytope $P(G)$, we would like to know
all the permutation groups $H$ such that $P(H)$ is combinatorially equivalent to $P(G)$. In the case of products 
the following result shows that we can reduce this question to each factor.

\begin{theorem}\label{product}
$P(G)$ is a combinatorial product of two polytopes $\Delta_1$ and $\Delta_2$ 
if and only if there are subgroups $H_1$ and $H_2$ in $G$ such that 
\begin{itemize}
\item[(a)] $P(H_i)$ is combinatorially equivalent to $\Delta_i$ for $i = 1,2$,
\item[(b)] $\su(H_1) \cap \su(H_2) = \emptyset$,
\item[(c)] $G = H_1 \times H_2$.
\end{itemize}
\end{theorem}

\begin{proof} 

The if-part is easy to see. We have to prove the only-if part.

Let $G \leq S_n$. By assumption, there is a map $v$ from the vertex set $\V(P(G)) = G$ to $\V(\Delta_1 \times \Delta_2) = 
\V(\Delta_1) \times \V(\Delta_2)$, 
inducing an isomorphism between the face lattices of $P(G)$ and of $\Delta := \Delta_1 \times \Delta_2$, which we also denote by $v$. Hence, any 
element $g \in G$ can be labeled as $v(g) = (v_1(g),v_2(g)) \in \V(\Delta)$ for unique vertices 
$v_1(g) \in \V(\Delta_1)$ and $v_2(g) \in \V(\Delta_2)$. We write $v(e) =: (e_1,e_2)$, and define 
$H_1 := \{g \in G \,:\, v_2(g) = e_2\}$, as well as $H_2 := \{g  \in G \,:\, v_1(g) = e_1\}$. 

We claim 
\begin{equation}
\su(H_1) \cap \su(H_2) = \emptyset.
\label{suppclaim}
\end{equation}

Let $h_1 \in H_1$ and $h_2 \in H_2$. We have $v(h_1) = (x_1,e_2)$ and $v(h_2) = (e_1,x_2)$ for $x_1 \in \V(\Delta_1)$ and 
$x_2 \in \V(\Delta_2)$. For $i=1,2$ let us denote by $F_i$ the smallest face of $\Delta_i$ containing $e_i$ and $x_i$. 
Let us define $g \in G$ with $v(g) = (x_1, x_2)$. 
Since $P(G)$ is via $v$ 
combinatorially equivalent to $\Delta$, the face $F_1 \times F_2 \prec \Delta$ is the smallest face of $\Delta$ 
containing $v(e)$ and $v(g)$.

By Corollary \ref{smallfacecoro} the face $F_g \prec P(G)$, satisfying $v(F_g) = F_1 \times F_2$, 
is centrally symmetric, and $h'_1 := g h_1^{-1} \in G$ is the antipodal vertex to $h_1$. 
Since $\su(h_1) \cap \su(h'_1) = \emptyset$, it suffices to show $h_2 = h'_1$.

Since $x_1 \times F_2$ is the smallest face of $\Delta$ containing 
$v(h_1)$ and $v(g)$, we get by central symmetry of $F_g$ that the smallest face of $\Delta$ 
containing $v(h'_1)$ and $v(e)$ has also $\card{\V(F_2)}$ vertices, thus 
$\card{\V(F_{h'_1})} = \card{\V(F_2)} = \card{\V(F_{h_2})}$. Note that by Theorem \ref{smallfacetheo} $h_1$, $h_2$ and $h'_1$ are subelements of $g$, 
thus determined by their support. 

Now, there are two cases, since $\su(g) = \su(h_1) \sqcup \su(h'_1)$ (here $\sqcup$ denotes the disjoint union).

\begin{enumerate}
\item $\su(h'_1) \subseteq \su(h_2)$: 

Then $h'_1$ is a subelement of $h_2$. However, 
$\card{\V(F_{h'_1})} = \card{\V(F_{h_2})}$ implies $F_{h'_1} = F_{h_2}$. Therefore, 
$h_2$ is also a subelement of $h'_1$, thus $h_2 = h'_1$, as desired.
\item $\su(h_1) \cup \su(h_2) \subsetneq \su(g)$: 

As in Remark \ref{reduction} for $S := \{e,h_1,h_2\}$ we define the matrix $M(S)$ and 
the face $F(S)$ of $P(G)$. Since $F_1 \times F_2$ is the smallest face of 
$\Delta$ containing $v(e),v(h_1),v(h_2)$, we get $F_g \subseteq F(S)$. However, by our assumption there exists some $i \in \su(g)$ with $i \not\in 
\su(h_1) \cup \su(h_2)$. Therefore, the only non-zero entry in the $i$th-row of $M(S)$ is on the diagonal, while the $i$th diagonal entry 
of the permutation matrix $g$ is zero. Hence, $\pro{M(S)}{g} < n$, a contradiction.
\end{enumerate}

This proves the claim (\ref{suppclaim}).
\smallskip

Hence, $\card{H_1 H_2} = \card{H_1} \card{H_2} = \card{\V(\Delta_1)} \card{\V(\Delta_2)} = \card{\V(P(G))} = \card{G}$. 
Therefore, $H_1 H_2 = G$. Moreover, this implies that $H_1$ consists precisely of all elements of 
$G$ that have disjoint support from all elements in $H_2$, hence is a subgroup. The analogous argument holds for $H_2$. 
Finally, $P(G)$, $P(H_1) \times P(H_2)$, and $\Delta_1 \times \Delta_2$ are combinatorially equivalent.
\end{proof}

As an application we classify those permutation groups whose $d$-dimen\-sional permutation polytopes 
have the maximal number of vertices.

\begin{corollary}
Let $d := \dim P(G)$. Then
\[\log_2 \card{G} \leq d \leq \card{G} - 1,\]
or equivalently
\[d + 1 \leq \card{G} \leq 2^d.\]
Moreover the following statements are equivalent:
\begin{enumerate}
\item $\card{G} = 2^d$,
\item $P(G)$ is combinatorially a $d$-cube,
\item $P(G)$ is lattice equivalent to $[0,1]^d$,
\item $G$ is effectively equivalent to $\<(1 2), \cdots, (2d\!-\!1 \;\: 2d)\> \leq S_{2 d}$.
\end{enumerate}
\label{maxtheo}
\end{corollary}

\begin{proof}

As was noted before, $P(G)$ is as a $d$-dimensional $0/1$-polytope that is combinatorially equivalent to a lattice subpolytope of $[0,1]^d$. 
Hence we get the inequalities. Moreover, it holds 
$(1) \lra (2)$, $(3) \ra (2)$ and $(4) \ra (2)$. From Theorem \ref{product} (in particular, statement (b)) we deduce $(2) \ra (3)$ and $(2) \ra (4)$.
\end{proof}

Here is another application of Theorem \ref{product}. For this note that by Theorem \ref{smallfacetheo} 
a permutation polytope $P(G)$ is simple if and only if there are dimension many indecomposable elements in $G$. 
Now, the main result of \cite{KW00} states that any simple $0/1$-polytope is a product of simplices. Therefore 
we can deduce from this geometric statement using Corollary \ref{simplex} a result in representation theory (since the dimension of $P(G)$ can 
be computed from the irreducible factors by Theorem \ref{dimension}).

\begin{corollary}
Let $\rho$ be a permutation representation of a group $G$. If $G$ contains precisely $\dim \,P(\rho)$ many indecomposable elements, then 
$G$ is the product of subgroups $H_1, \ldots, H_l$ with mutually disjoint support, where $\rho$ restricted to any $H_i$ ($i = 1, \ldots, l$) 
is stably equivalent to the regular representation of $H_i$.
\end{corollary}

\subsection{Pyramids}

In experiments one observes that most faces of permutation polytopes are actually pyramids over lower dimensional faces. Here, we prove 
that for any face $F$ of a permutation polytope there exists a permutation polytope having a face that is combinatorially a pyramid over $F$. 

\begin{theorem}
Let $G \leq S_n$. Then there is a permutation group $E \leq S_{2n}$, with $E \cong G \times \Z_2$, 
such that there is a face of $P(E)$ which is combinatorially a pyramid over $P(G)$.
\label{pyramid}
\end{theorem}

\begin{proof}

Let $G \leq S_n$. Embedding the product of permutation groups $S_n \times S_n$ into $S_{2n}$, 
we define $H := \{(\sigma, \sigma) \,:\, \sigma \in G\} \leq S_{2n}$. 
Then $H$ is a subgroup of $S_{2n}$, and effectively equivalent 
to $G$. We define an involution 
\[p := (1 \;\;n\!\!+\!\!1)\;\; (2 \;\;n\!\!+\!\!2) \; \cdots \; (n\;\; 2n) \in S_{2n}.\]
Then $p$ commutes with each element in $H$, moreover $H \cap \<p\> = \{e\}$. 
Hence, $E := H \<p\>$ is a subgroup of $S_{2n}$ and isomorphic to $H \times \<p\>$. 
For $S := H \sqcup \{p\} \subseteq E$ we define as in Remark \ref{reduction} 
the $2n \times 2n$-matrix $M(S)$ defining a face $F(S)$ of $P(E)$. 
We claim that $F(S)$ is a pyramid over $\conv(H)$. 

First let us show that $\V(F(S)) = S$. Assume that there is some $h p \in F(S)$ with $h \not= e$. Let $h = (\sigma, \sigma)$ for 
$\sigma \in G$, $\sigma \not= e$. Assume $\sigma$ maps $1$ to $2$. 
Then $h p$ maps $1$ to $n+2$. However, this implies that $\pro{M(S)}{h p} < 2n$, a contradiction.

Now, it remains to show that $H$ is the set of vertices of a face of $F(S)$. 
As in Remark \ref{reduction} we define the face $F(H)$ of $P(E)$. By construction $p \not\in F(H)$. 
Then $F(H) \cap F(S)$ is a face of $F(S)$ which contains $H$ but not $p$.
\end{proof}

\begin{corollary}
Pyramids over faces of permutation polytopes appear as faces of permutation polytopes.
\label{pyracor}
\end{corollary}

\subsection{Free sums}

Recall that free sums are the combinatorially dual operation to products. 
For instance the free sum of $d$ intervals is a $d$-crosspolytope, i.e., 
the centrally symmetric $d$-polytope with the minimal number $2 d$ of vertices. 

In general we cannot expect that free sums of arbitrary faces of permutation polytopes are again faces of permutation polytopes. 
Corollary \ref{smallfacecoro} shows that already faces that are bipyramids have to be necessarily centrally symmetric. 
However, we can explicitly construct the following centrally symmetric polytopes as faces.

\begin{theorem}
Let $l,d$ be natural numbers. There exists a face of a permutation polytope that is combinatorially 
the free sum of an $l$-crosspolytope and a $d$-cube.
\label{constr}
\end{theorem}

\begin{proof}

Since for $l \geq 1$ an $l$-crosspolytope is the free sum of an $(l-1)$-crosspolytope and a $1$-cube, we may assume $l,d \geq 1$. 
We set $n_0 := 3 d$. Let $z_1, \ldots, z_d \in S_{n_0}$ be disjoint $3$-cycles. We define $G_0 := \<z_1, \ldots, z_d\>$, an elementary abelian $3$-group 
of order $3^d$. It contains $g_0 := z_1 \cdots z_d$. Let $V_0 := \{z_1^{k_1} \cdots z_d^{k_d} \;:\; k_i \in \{0,1\}\}$. 
This is precisely the set of subelements of $g_0$ in $G_0$, hence, by Theorem \ref{smallfacetheo} 
$V_0$ is the vertex set of a face $F_0$ of $P(G_0)$. 
This face $F_0$ is combinatorially a $d$-cube. 

Now, we proceed by induction for $i= 1,\ldots, l$. We define 
\[n_i := 2 n_{i-1},\quad H_i := \{(g,g) \;:\; g \in G_{i-1}\} \leq S_{n_i}, \quad g_i := (g_{i-1}, g_{i-1}) \in H_i,\] 
\[p_i := (g_{i-1}, e) \in S_{n_i}, \quad p'_i := (e, g_{i-1}) \in S_{n_i}, \quad G_i := H_i \<p_i\> \leq S_{n_i}.\]
Note that $H_i, G_i$ are elementary abelian $3$-groups. Moreover, $G_i \cong H_i \times \<p_i\>$. 
Let $V_i$ be the set of subelements of 
$g_i$ in $G_i$, and $F_i$ the smallest face of $P(G_i)$ containing $e$ and $g_i$, thus $V_i$ is the vertex set of $F_i$ by Theorem \ref{smallfacetheo}. 
By induction hypothesis we know that $F_{i-1}$ is combinatorially the free sum of an $i-1$-crosspolytope and a $d$-cube. 
We show that $F_i$ is a bipyramid over $F_{i-1}$ with apexes $\{p_i, p'_i\}$. We claim
\begin{equation}
V_i = \{(v,v) \;:\; v \in V_{i-1}\} \sqcup \{p_i, p'_i\}.
\label{vclaim}
\end{equation}

Let $(v,v') \in V_i$. Then $v,v'$ are subelements of $g_{i-1}$ in $G_{i-1}$. Assume $v \not= v'$, in particular $(v,v') \not\in H_i$. 
Since $(v,v') \in G_i$ and $p_i$ has order three, 
we have either (i) $(v,v') = h p_i$ or (ii) $(v,v') = h p^{-1}_i$ for some $h \in H_i$. Let $h = (g,g)$ for $g \in G_{i-1}$. 
This yields either $(v,v') = (g g_{i-1}, g)$ or $(v,v') = (g g^{-1}_{i-1}, g)$, in particular, $g=v'$ is a subelement of $g_{i-1}$. 
Now, the case (i) implies $g=e$, thus $(v,v') = p_i$, since otherwise $g g_{i-1}$ cannot be a subelement of $g_{i-1}$, 
because $g_{i-1}$ decomposes into disjoint $3$-cycles by construction. 
In the case (ii) we see analogously $g=g_{i-1}$, thus $(v,v') = p'_i$. This proves the claim (\ref{vclaim}).
\smallskip

Since by Corollary \ref{smallfacecoro} 
$F_i$ is centrally symmetric with antipodal vertices $p_i$ and $p'_i$, and since $p_i, p'_i \not\in \aff H_i$ by construction, 
the claim implies that the face $F_i$ is a bipyramid over the convex hull of all elements 
$(v,v)$ (for $v \in V_{i-1}$), which is affinely equivalent to $F_{i-1}$.
\end{proof}

\begin{corollary}
For any dimension $d$ there is a face of a permutation polytope that is a $d$-crosspolytope.
\label{crosscoro}
\end{corollary}

\begin{example}{\rm
As an illustration of the proof of Theorem \ref{constr}, we show 
how to obtain the octahedron as a face of a permutation polytope. The octahedron is the free sum 
of an interval and a square, so $l=1$ and $d=2$. Therefore, we define $z_1 := (1 2 3), z_2 := (4 5 6)$, and 
$G_0 := \<z_1, z_2\> \leq S_6$. Now, we set $z'_1 := (7 8 9), z'_2 := (10\, 11\, 12)$, and 
$G_1 := \< z_1 z'_1, z_2 z'_2, z_1^{-1} z_2 z'_1\> \leq S_{12}$. Then for $g_1 := z_1 z_2 z'_1 z'_2$ 
the face $F_1 := F_{g_1}$ of $P(G)$ is an octahedron with 
the vertex set $V_1 = \{e, z_1 z'_1, z_2 z'_2, z_1 z_2 z'_1 z'_2, z_1 z_2, z'_1 z'_2\}$.
\label{example-octahedron}
}
\end{example}

\subsection{Centrally symmetric permutation polytopes}

We will establish a one-to-one correspondence between centrally
symmetric permutation polytopes on one hand and certain subspaces of
$\F^r$ on the other. We will liberally identify sets $I \subseteq [r] := \{1,
\ldots, r\}$ with their incidence vectors $I \in \F^r$.

Suppose $P(G)$ is centrally symmetric, and let $g_0 = z_1 \circ
\cdots \circ z_r \in G$ be the vertex opposite to $e$. Then $P(G) =
F_{g_0}$, and by Theorem~\ref{smallfacetheo}, every element of $G$
is a subelement of $g_0$. Applied to $g_0^{-1} = z_1^{-1} \circ
\cdots \circ z_r^{-1}$ this yields that for $i = 1, \ldots, r$ we have
$z_i^{-1} = z_i$, so that $z_i$ is a cycle of length two. Hence, the
elements of $G$ have order two. They can be written as $z_I := \Pi_{i
  \in I} z_i$ for certain $I \subseteq [r]$. (For example,
$z_\emptyset = e$, and $z_{[r]} = g_0$.)

Multiplication in $G$ corresponds to addition (symmetric difference)
in $\F^r$. This means that the set of $I$ such that $z_I \in G$ is
a subspace of $\F^r$. Conversely, given such a subspace containing
the all-ones-vector $[r]$, we obtain a centrally symmetric permutation
polytope. We have just proved the second part of the following
proposition. The first part follows from
Corollary~\ref{smallfacecoro}.

\begin{prop}
  $P(G)$ is centrally symmetric if and only if there is a pair of
  vertices such that $P(G)$ is the smallest face containing them.

  In this case, $G$ is an elementary abelian $2$-group, in particular
  the number of vertices of $P(G)$ is a power of two. \label{twoprop}
\end{prop}

If $G$ is an elementary abelian $2$-group, $P(G)$ does not have to be
centrally symmetric, see Theorem \ref{permpolytheo}.

Continuing our discussion of centrally symmetric permutation
polytopes, we can reorder the elements of $[n]$ ($G \subseteq S_n$) so
that the matrices in $G$ are block diagonal with blocks
$\left(
  \begin{smallmatrix}
    1&0\\0&1
  \end{smallmatrix}
\right)$
or
$\left(
  \begin{smallmatrix}
    0&1\\1&0
  \end{smallmatrix}
\right)$ (plus $n-2r$ blocks $(1)$).
Projecting to the $r$ upper right entries of the $2 \times 2$ blocks
is a lattice isomorphism to $\Z^r$. Reduction mod $2$ yields, again,
our subspace.

\begin{prop}
  If $P(G)$ is centrally symmetric, then the free sum of $P(G)$ with itself 
  is again a permutation polytope.\label{directsumprop}
\end{prop}

\begin{proof}
  Suppose $P(G)$ is centrally symmetric with corresponding subspace $V
  \subseteq \F^r$. Define a new subspace $$\widetilde V \ := \ \{
  (I,I) \ : \ I \in V \} \cup \{ (I, [r] - I) \ : \ I \in V \} \
  \subset \ \F^r \times \F^r$$
  We claim that the permutation polytope $P(\widetilde G)$ of the 
  corresponding permutation group $\widetilde G$ realizes the free sum 
  of $P(G)$ with itself.

  We work with the $2r$ upper right entries $(x,y) \in \R^r \times
  \R^r$. Both, the diagonal embedding $x \mapsto (x,x)$ and the
  ``anti-diagonal'' embedding $x \mapsto (x,1-x)$ include $P(G)$ into
  $P(\widetilde G)$ and exhaust all vertices. Their images intersect
  in the unique common interior point $\frac{1}{2}(1,1)$.
\end{proof}

While in any dimension $d$ there exists a permutation polytope which
is a $d$-cube, this is not true for its combinatorial dual, the
$d$-dimensional crosspolytope.

\begin{theorem}
  There  is a  $d$-dimensional  permutation  polytope  $P(G)$ that  is
  combinatorially  a cross\-po\-ly\-to\-pe if   and only if  $d$ is  a
  power of two.  In this case,  the effective equivalence class of $G$
  is uniquely determined.
  \label{cross}
\end{theorem}

\begin{proof}
  The fact that $d$ must be a power of two follows from
  Proposition~\ref{twoprop}. In that case, existence follows from
  Proposition~\ref{directsumprop}. So we only need to show uniqueness.

  Let $G$ be a permutation group such that $P(G)$ is a $d$-dimensional
  cross\-poly\-tope, $d=2^k$. Then $|G| = 2^{k+1}$, and our subspace $V
  \subset \F^r$ has dimension $k+1$. 
  Choose generators $g_0, g_1, \ldots, g_k$ of $G$, i.e., a
  basis of $V$.
  Consider the $k \times r$ matrix with rows $g_1, \ldots, g_k \in
  \F^r$. If there are two equal columns $i,j$ then we can omit $z_i$
  without changing the effective equivalence class. On the other hand,
  there can be at most $2^k$ different columns.
  As remarked above, we can embed $P(G) \hookrightarrow \R^r$. So, in
  particular, $r \ge d = 2^k$. Hence, $r=d$, and this matrix simply
  lists all possible $0/1$ vectors. It is, up to permutation of the
  columns, uniquely defined.
\end{proof}

\begin{example}{\rm 
From the proof we get an explicit description of the permutation groups 
defining $d$-dimensional crosspolytopes. For instance, let $d=4$. 
Since $r=d=4$, we have $g_0 = z_1 \circ z_2 \circ z_3 \circ z_4$, 
where we may choose $z_1 = (1 2)$, $z_2 = (3 4)$, $z_3 = (5 6)$, $z_4 = (7 8)$. 
Since $k=2$, we get $G = \< g_0, g_1, g_2\>$, where $g_1$ and $g_2$ are given 
by the rows of the following matrix consisting of all possible vectors in 
$\{0,1\}^2$:
\[\begin{pmatrix} 
0 & 1 & 0 & 1\\
0 & 0 & 1 & 1
\end{pmatrix}\]
So, $g_1 = z_2 z_4$, $g_2 = z_3 z_4$, and $G = \< (1 2) (3 4) (5 6) (7 8), (3 4) (7 8), (5 6) (7 8)\>$.
\label{example-4-crosspolytope}
}
\end{example}

\section{Classification results in low dimensions} \label{sec:classification}

\subsection{Classification of $\leq 4$-dimensional permutation polytopes}

We would like to classify all permutation polytopes of given small dimension $d$. 
For this we take a look at the list of Aichholzer \cite{Aic07} of 
combinatorial types of $0/1$-polytopes in small dimension. For any such $0/1$-polytope we first check whether it has constant vertex degree 
and satisfies the condition of Corollary \ref{smallfacecoro}. Then we go through the list 
of groups of size equal to the given number of vertices. Now, using the theoretical results of the previous section we can deduce 
from the combinatorial structure of the polytope whether this polytope can be realized as a permutation polytope, and 
even determine all respective permutation groups up to effective equivalence.

\begin{theorem}
Table \ref{pp-table1} contains the list of all permutation groups 
$G$ with $d$-dimensional permutation polytope $P(G)$ for $d \leq 4$ up to effective equivalence.
\label{permpolytheo}
\end{theorem}
\begin{table}[t]
\setlength{\belowcaptionskip}{5pt}
\renewcommand{\arraystretch}{1.2}
\begin{tabular}[t]{l|l|l}
Combin.\ type of $P(G)$ & Isom.\ type of $G$ & Effective equiv.\ type of $G$\\\hline\hline
triangle & $\Z/3\Z$ & $\<(1 2 3)\>$\\
square & $(\Z/2\Z)^2$ & $\<(1 2), (3 4)\>$\\
\cline{1-3}tetrahedron & $\Z/4\Z$ & $\<(1 2 3 4)\>$\\
tetrahedron & $(\Z/2\Z)^2$ & $\<(1 2)(3 4), (13)(24)\>$\\
triangular prism & $\Z/6\Z$ & $\<(1 2), (3 4 5)\>$\\
cube & $(\Z/2\Z)^3$ & $\<(1 2), (3 4), (5 6)\>$\\
\cline{1-3}$4$-simplex & $\Z/5\Z$ & $\<( 1 2 3 4 5)\>$\\
$B_3$ & $S_3$ & $\<(1 2), (1 2 3)\>$\\
prism over tetrahedron & $\Z/2\Z \times \Z/4\Z$ & $\<(1 2 3 4),(5 6)\>$\\
prism over tetrahedron & $(\Z/2\Z)^3$ & $\<(1 2) (3 4), (1 3) (24), (5 6)\>$\\
$4$-crosspolytope & $(\Z/2\Z)^3$ & $\<(1 2) (3 4), (3 4) (7 8), (5 6) (7 8)\>$\\
product of triangles & $(\Z/3\Z)^2$& $\<(1 2 3), (4 5 6)\>$\\
prism over triang.\ prism & $\Z/6\Z \times \Z/2\Z$ & $\<(1 2), (3 4 5), (6 7)\>$\\
$4$-cube & $(\Z/2\Z)^4$ & $\<(1 2), (3 4), (5 6), (7 8)\>$
\end{tabular}
\caption{\label{pp-table1}Permutation polytopes in dimension $\leq 4$}
\end{table}
\begin{proof}
$d=2$: The triangle and the square are the only two-dimensional $0/1$-polytopes. 
If $P(G)$ is a triangle, then use Corollary \ref{simplex}. 
If $P(G)$ is a square, then use Corollary \ref{maxtheo}.

$d=3$: There are $4$ combinatorial types of three-dimensional $0/1$-polytopes 
with constant vertex degree satisfying the condition of Corollary \ref{smallfacecoro}.

\begin{enumerate}
\item $P(G)$ is a tetrahedron: Then use Corollary \ref{simplex}.

\item $P(G)$ is a triangular prism: Then use Theorem \ref{product}.

\item $P(G)$ is an octahedron: Then $\card{G} = 6$, however by Proposition \ref{twoprop} $\card{G}$ has to be a power of two, a contradiction.

\item $P(G)$ is a cube: Then use Corollary \ref{maxtheo}.
\end{enumerate}

$d=4$:  There are $9$ combinatorial types of four-dimensional $0/1$-polytopes 
with constant vertex degree satisfying the condition of Corollary \ref{smallfacecoro}.

\begin{enumerate}
\item $P(G)$ is a $4$-simplex: Then use Corollary \ref{simplex}. 

\item $P(G)$ is combinatorially equivalent to the Birkhoff polytope $B_3 = P(S_3)$. 
Then $|G| = 6$. There are two cases:

If $G \cong \Z/6\Z$, then either $G$ is generated by an element of order $6$, thus $G$ is effectively equivalent to $\< (1 2 3 4 5 6) \> \leq S_6$, 
or $G$ is generated by two elements of orders $2$ and $3$. In the latter case, the two elements necessarily have disjoint support, 
since $G$ is abelian. In the first case $P(G)$ is a $5$-simplex, in the second case $P(G)$ is a prism. Both cases yield contradictions.

Hence, $G \cong S_3$. Now, $\Irr(S_3) = \{1_{S_3}, \rho_1, \rho_2\}$ with $\deg \rho_1 = 1$ and $\deg \rho_2 = 2$. 
By the dimension formula the permutation representation $\rho$ associated to the permutation group $G$ can only have $\rho_2$ as an 
irreducible factor. Therefore, Theorem \ref{thm:stablyEquivalence} implies that $G$ is effectively equivalent to $\<(1 2), (1 2 3)\>$.

\item $P(G)$ is a prism over a tetrahedron: Use Theorem \ref{product} and the classification for $d=3$.

\item $P(G)$ is a $4$-crosspolytope: See Example \ref{example-4-crosspolytope}.

\item $P(G)$ is a product of two triangles: Use Theorem \ref{product}.

\item $P(G)$ is a prism over the triangular prism: Use Theorem \ref{product} and the classification for $d=3$.

\item $P(G)$ is a $4$-cube: Use Corollary \ref{maxtheo}.

\item $P(G)$ is a prism over the octahedron: Then 
Proposition \ref{twoprop} yields that the number of vertices has to be a power of two, but $P(G)$ has $12$ vertices, a contradiction.

\item $P(G)$ is a hypersimplex: 
Any vertex is contained in precisely three facets that are octahedra. Since the inversion map on $G$ is given 
by the transposition map on $P(G)$, it induces an automorphism of $P(G)$ of order two, and hence, since $3$ is odd, 
there has to be an octahedron $F$ 
that contains $e$ and whose vertex set is invariant under inversion. 

Let $g$ be the unique vertex of $F$ opposite to $e$. Hence, $g$ is fixed by the inversion map, so 
$g=g^{-1}$. As already noted in the proof of Proposition \ref{twoprop}, 
this yields that any vertex of $F=F_g$ (besides $e$) has order two, so there are at least five elements of order two in $G$. 
On the other hand, since $5$ divides $\card{G} = 10$, there exists a subgroup of order $5$, so we conclude that 
there are precisely four elements of $G$ that have order~$5$. 

Now take $F'$ ($\not= F$) as one of the other two octahedra that contain $e$. We denote by $h$ the unique vertex $h$ ($\not= g$) in $F'$ opposite to $e$. 
As just seen, the order of $h$ has to be $5$. Since $h$ has to be a product of disjoint cycles of order $5$, any 
subelement of $h$ ($\not= e$) also has order $5$, so by Theorem
\ref{smallfacetheo} there are at least five vertices of $F'=F_h$ of
order $5$, a contradiction.
\qedhere
\end{enumerate}
\end{proof}

\begin{remark}{\rm Another approach following 
Theorem \ref{thm:stablyEquivalence} and Corollary \ref{maxtheo} would be to determine all abstract groups $G$ of order $\leq 2^d$ and to 
calculate the finite set $\Irr(G)$. Then for any subset $S \subseteq \Irr(G)$ 
it would be enough to find, if possible, some permutation representation $\rho$ with $\Irr(\rho) = S$. However, 
the last task seems to be neither practically nor theoretically easy to achieve, compare \cite{BP98}. 
}
\end{remark}

\subsection{Classification of $\leq 4$-dimensional faces}

Compared to the classification of permutation polytopes the question whether a given $0/1$-polytope is combinatorially equivalent to the 
face of some permutation polytope is much more difficult. If the answer is supposed to be positive, 
then one has to construct an explicit permutation group, the dimension of whose permutation polytope might increase dramatically. 
A systematic way to perform this task is yet to be discovered. However, to show that the answer is negative is even more challenging, since we 
lack good combinatorial obstructions of the type given in Corollary \ref{smallfacecoro}. 

\begin{theorem}
The following list contains all combinatorial types of $d$-dimensional $0/1$-polytopes for $d \leq 4$ 
that may possibly appear as faces $F$ of some permutation polytope $P(G)$:

\begin{enumerate}
\item[$d=2$:] There are $2$ combinatorial types realized as $F$: triangle and square. Both appear as faces of Birkhoff polytopes.
\item[$d=3$:] There are $5$ combinatorial types realized as $F$: tetrahedron, square pyramid, triangular prism, cube, and 
octahedron. The first four appear as faces of Birkhoff polytopes.
\item[$d=4$:] There are $21$ combinatorial types that may possibly appear as $F$:
\begin{enumerate}
\item $11$ of these appear as faces of Birkhoff polytopes: $4$-simplex, pyramid over square pyramid, Birkhoff polytope $B_3$ (free sum of 
two triangles), pyramid over prism 
over triangle, wedge $W$ over base edge of square pyramid, pyramid over cube, prism over tetrahedron, 
product of two triangles, prism over square pyramid, product of triangle and square, $4$-cube.
\item $8$ of these can be realized as $F$: $4$-crosspolytope, prism over octahedron, pyramid over octahedron, bipyramid over cube, wedge 
over the facet of an octahedron, dual of $W$ (see (a)), 
hypersimplex (the combinatorial type of $\{x \in [0,1]^5 \,:\, \sum_{i=1}^5 x_i = 2\}$), 
and one special $0/1$-polytope $P$.
\item $2$ of these are given by special $0/1$-polytopes $Q_1,Q_2$, where it is unknown, if they have a realization as $F$.
\end{enumerate}
\smallskip

The description of the combinatorial types of $P$, $Q_1$, $Q_2$ can be found in Table \ref{pp-table2}.

\begin{table}[b]
  \small
  \centering
\setlength{\abovecaptionskip}{5pt}
\setlength{\belowcaptionskip}{3pt}
\setlength{\subfigcapskip}{5pt}
  \caption{\label{pp-table2}Vertex-Facet-Incidences of $P$, $Q_1$, $Q_2$}	
\label{tab:vifs-prod} \subtable[$P$: f-vector $(8,21,22,9)$]{
    \centering
    \begin{tabular*}{.3\textwidth}{@{\extracolsep{-5pt}}
        >{$}r<{$}>{$}r<{$}>{$}r<{$}>{$}r<{$}>{$}r<{$}>{$}r<{$}>{$}r<{$}}
      [&0& 2& 5&\multicolumn{2}{>{$}l<{$}}{6\;]} \\[.05cm]
      [&0& 2& 4&\multicolumn{2}{>{$}l<{$}}{6\;]}\\[.05cm]
      [&0& 1& 2& 3&\multicolumn{2}{>{$}l<{$}}{4\;]}\\[.05cm]
      [&0& 1& 2& 3&\multicolumn{2}{>{$}l<{$}}{5\;]}\\[.05cm]
      [&1& 3& 4&\multicolumn{2}{>{$}l<{$}}{7\;]}\\[.05cm]
      [&1& 3& 5&\multicolumn{2}{>{$}l<{$}}{7\;]}\\[.05cm]
      [&2& 3& 4& 6&\multicolumn{2}{>{$}l<{$}}{7\;]}\\[.05cm]
      [&2& 3& 5& 6&\multicolumn{2}{>{$}l<{$}}{7\;]}\\[.05cm]
      [&0& 1& 4& 5& 6& 7\;]
    \end{tabular*}
  }
  \subtable[$Q_1$: f-vector $(7,19,23,11)$]{
    \centering
    \begin{tabular*}{.3\textwidth}{@{\extracolsep{-5pt}}
        >{$}r<{$}>{$}r<{$}>{$}r<{$}>{$}r<{$}>{$}r<{$}>{$}r<{$}}
      [&0 &1& 2& 3 & 4\;]\\[.05cm]
      [&1 &2& 3& 4 & 5 \;]\\[.05cm]
      [&0 &1& 2&\multicolumn{2}{>{$}l<{$}}{5\;]}\\[.05cm]
      [&1 &3& 5&\multicolumn{2}{>{$}l<{$}}{6\;]}\\[.05cm]
      [&0 &1& 3&\multicolumn{2}{>{$}l<{$}}{6\;]}\\[.05cm]
      [&0 &1& 5&\multicolumn{2}{>{$}l<{$}}{6\;]}\\[.05cm]
      [&0 &2& 5&\multicolumn{2}{>{$}l<{$}}{6\;]}\\[.05cm]
      [&2 &4& 5&\multicolumn{2}{>{$}l<{$}}{6\;]}\\[.05cm]
      [&0 &2& 4&\multicolumn{2}{>{$}l<{$}}{6\;]}\\[.05cm]
      [&3 &4& 5&\multicolumn{2}{>{$}l<{$}}{6\;]}\\[.05cm]
      [&0 &3& 4&\multicolumn{2}{>{$}l<{$}}{6\;]}
    \end{tabular*}
  }
  \subtable[$Q_2$: f-vector $(8,25,32,15)$]{
    \centering
    \begin{tabular*}{.3\textwidth}{@{\extracolsep{-5pt}}
        >{$}r<{$}>{$}r<{$}>{$}r<{$}>{$}r<{$}>{$}r<{$}>{$}r<{$}>{$}r<{$}}
      [&0 &1& 2& 3 & 4 & 5\;]\\[.05cm]
      [&2 &4& 5&\multicolumn{2}{>{$}l<{$}}{6\;]}\\[.05cm]
      [&1 &3& 5&\multicolumn{2}{>{$}l<{$}}{6\;]}\\[.05cm]
      [&1 &2& 5&\multicolumn{2}{>{$}l<{$}}{6\;]}\\[.05cm]
      [&0 &1& 2&\multicolumn{2}{>{$}l<{$}}{6\;]}\\[.05cm]
      [&0 &1& 3&\multicolumn{2}{>{$}l<{$}}{7\;]}\\[.05cm]
      [&0 &1& 6&\multicolumn{2}{>{$}l<{$}}{7\;]}\\[.05cm]
      [&0 &2& 6&\multicolumn{2}{>{$}l<{$}}{7\;]}\\[.05cm]
      [&0 &2& 4&\multicolumn{2}{>{$}l<{$}}{7\;]}\\[.05cm]
      [&0 &3& 4&\multicolumn{2}{>{$}l<{$}}{7\;]}\\[.05cm]
      [&2 &4& 6&\multicolumn{2}{>{$}l<{$}}{7\;]}\\[.05cm]
      [&1 &3& 6&\multicolumn{2}{>{$}l<{$}}{7\;]}\\[.05cm]
      [&3 &5& 6&\multicolumn{2}{>{$}l<{$}}{7\;]}\\[.05cm]
      [&4 &5& 6&\multicolumn{2}{>{$}l<{$}}{7\;]}\\[.05cm]
      [&3 &4& 5&\multicolumn{2}{>{$}l<{$}}{7\;]}
    \end{tabular*}
  }
\end{table}

\end{enumerate}
\label{facestheo}
\end{theorem}

\begin{proof}

By \cite{BS96} any 
$d$-dimensional face $F$ of some Birkhoff polytope is already realized in $B_{2 d}$. Hence, by looking at the faces of $B_6,B_8$, we 
find all combinatorial types of $\leq 4$-dimensional faces of Birkhoff polytopes. 

Let us now consider the general case of a $d$-dimensional 
face of a permutation polytope.

$d=2$: The triangle and the square are the only two-dimensional $0/1$-polytopes.

$d=3$: There are $5$ combinatorial types of three-dimensional $0/1$-polytopes 
that satisfy the condition of Corollary \ref{smallfacecoro}. 
Since the first four are realized as Birkhoff polytopes, 
we only have to deal with the octahedron. This was done in Example \ref{example-octahedron}.

$d=4$: There are $21$ combinatorial types of four-dimensional $0/1$-polytopes 
that satisfy the condition of Corollary \ref{smallfacecoro}. Of these, $11$ can be realized as faces of 
Birkhoff polytopes. Here are the remaining $10$ cases:

\begin{enumerate}
\item The $4$-crosspolytope: See Table \ref{pp-table1}.
\item The prism over an octahedron: The octahedron is a face of a permutation polytope, so also the prism is.
\item Pyramid over octahedron: See Corollary \ref{pyracor} and the classification for $d=3$.
\item Bipyramid over cube: See Theorem \ref{constr}.
\item The dual of $W$: 
Let $a_1,a_2,a_3,a_4,b_1,b_2,b_3,b_4$ be eight $3$-cycles with pairwise disjoint support, realized as elements in $S_{24}$. 
Let $e_1$ be an involution that exchanges $a_1$ and $a_2$, i.e., for $a_1 = (1 2 3)$ and $a_2 = (4 5 6)$, we 
define $e_1 := (1 4) (2 5) (3 6)$. In the same way we define $e_2$ as the involution exchanging $b_1$ and $b_2$, and 
$d_1$ exchanging $a_3$ and $a_4$, and $d_2$ exchanging $b_3$ and $b_4$. Then we define 
$v_1 := a_1 a_2 a_3 a_4$, $v_2 := b_1 b_2 b_3 b_4$, $v_3 := d_1 d_2$, and $v_4 := e_1 e_2$. 
Let $G := \<v_1,v_2,v_3,v_4\> \leq S_{24}$. Then $\card{G} = 36$. As in Remark \ref{reduction} we define 
the face $F := F(\{e,v_1,v_2,v_3,v_4\})$ of $P(G)$. 
Now, we check using \texttt{GAP} and \texttt{polymake} that the combinatorial type of $F$ is indeed the dual of $W$.
\item $P$: 
Let $a_1,a_2,b_1,b_2,c_1,c_2$ be six $3$-cycles with pairwise disjoint support, realized as elements in $S_{18}$. 
Let $e_1$ be an involution that exchanges $b_1$ and $b_2$, as before, 
and $e_2$ an involution that exchanges $c_1$ and $c_2$. Then we define $v_1 := a_1 a_2$, $v_2 := b_1 b_2 c_1 c_2$, $v_3 := a_1 b_1 b_2$, $v_4 := e_1 e_2$. 
Let $G := \<v_1,v_2,v_3,v_4\> \leq S_{18}$. Then $\card{G} = 54$. As in Remark \ref{reduction} we define 
the face $F := F(\{e,v_1,v_2,v_4\})$ of $P(G)$. Now, we check that the combinatorial type of $F$ is indeed $P$.
\item The wedge over the facet of an octahedron: 
Let $a_1,a_2$,$b_1,b_2$,$c_1,c_2$,$e_1$,\\$v_1,v_2,v_3$ be defined as for $P$. 
However $v_4 := e_1$. Then $G:= \<v_1,v_2,v_3,v_4\> \leq S_{18}$ with $\card{G} = 54$. We define $F$ as before. 
Now, we check that the combinatorial type of $F$ is indeed as desired.
\item Hypersimplex: Let $a_1,a_2,a_3,a_4,a_5$ be five $3$-cycles with pairwise disjoint support, realized as elements in $S_{15}$. 
We define $v_1 := a_1 a_2$, $v_2 := a_3 a_4$, $v_3 := a_1 a_3$, $v_4 := a_1 a_5$. Let $G:= \<v_1,v_2,v_3,v_4\> \leq S_{15}$. Then 
$\card{G} = 81$. As in Remark \ref{reduction} we define the face $F := F(\{e,v_1,v_2,v_4\})$ of $P(G)$. 
Then $F$ has dimension five, and we check that $F$ is a pyramid over the hypersimplex.
\item $Q_1$: 
We could not find a permutation group $G$ with a face $F$ of $P(G)$ whose combinatorial type coincides with the one of $Q_1$.
\item $Q_2$: 
As for $Q_1$.\qedhere
\end{enumerate}
\end{proof}

\begin{conjecture}
{\rm There are no four-dimensional faces of permutation polytopes having the combinatorial type of $Q_1$ or $Q_2$ (see Table 
\ref{pp-table2}).
}
\end{conjecture}

\begin{remark}{\rm 
The combinatorics of faces of permutation polytopes is in general much more complex than the one 
of faces of Birkhoff polytopes. For instance, any facet of a Birkhoff polytope is given as the set of matrices with 
entry $0$ (resp. $1$) at a fixed position $i,j$. Hence, any face of a Birkhoff polytope has the strong property that 
the vertices in the complement of a facet form a face. 
In dimension $d \leq 4$ all of the polytopes in Theorem \ref{facestheo} also possess this property, except the 
special $0/1$-polytopes $Q_1$, $Q_2$, which stresses the exceptional role of these polytopes. However, one should not 
jump to the wrong conclusion that this might be a necessary condition on a polytope to be a face of a permutation polytope. 
In dimension $d \geq 9$ there are examples 
of permutation polytopes that have facets whose complement is not even a subset of a proper face, see \cite{BHNP07}.
\label{strangefaces}
}\end{remark}

\section{Open questions and conjectures} \label{sec:questions}

\subsection{Permutation polytopes}

Inspired by an embedding result on Birk\-hoff faces in \cite{BS96} we propose the following 
daring conjecture (in a weak and strong version), where the bound $2d$ would be sharp 
as the example of the $d$-cube shows, see Corollary \ref{maxtheo}. The existence of some bound follows from Proposition \ref{01bound}.

\begin{conjecture}{\rm 
    Let $P$ be a $d$-dimensional permutation polytope. Then there
    exists a permutation group $G \leq S_{2 d}$ such that $P(G)$ is
    combinatorially equivalent (or stronger, lattice equivalent) to $P$. 
\label{combconjecture}
  }
\end{conjecture}

An even more natural formulation is given be the next conjecture, which was checked for $d \leq 4$ using Theorem \ref{permpolytheo}. 
The statement may be phrased purely in terms of representation theory thanks to the dimension formula, see Theorem \ref{dimension}.

\begin{conjecture}{\rm 
Let $\rho$ be a permutation representation of a finite group $G$ with $d := \dim P(\rho)$. 
Then there exists a stably equivalent permutation representation $\rho' \colon G \to S_n$ such that $n \leq 2 d$.
}
\end{conjecture}

The truth of this statement would imply the weak part of Conjecture \ref{combconjecture}. 

\smallskip

The Birkhoff polytope $B_3$ is given by the full symmetric group $S_3$. From 
Theorem \ref{permpolytheo} we observe that $S_3$ is essentially the only permutation group having $B_3$ as its 
permutation polytope.

\begin{conjecture}{\rm 
    Let $P(G)$ be a permutation polytope such that $P(G)$ is
    combinatorially equivalent to the Birkhoff polytope $B_n$ for some
    $n$. Then the permutation group $G$ is effectively equivalent to
    $S_n$. 
  }
\end{conjecture}

Any element of a permutation group $G$ induces by left, respectively, by right multiplication an affine automorphism 
of $P(G)$. If $G$ is not abelian, this implies that there are more affine automorphisms of $P(G)$ than elements of $G$. 
We conjecture this to be true also in the abelian case, except if $\card{G} \leq 2$.

\begin{conjecture}{\rm Let $G$ be an abelian permutation group of order $\card{G} > 2$. 
Then the group of affine automorphisms of $P(G)$ contains more elements than $G$.
}
\end{conjecture}

\subsection{Faces of permutation polytopes}

Observing the structure of centrally symmetric faces for $d \leq 4$ gives rise to the following question.

\begin{question}{\rm 
Is there a centrally symmetric face of a permutation polytope that is not composable as products or free sums 
of lower dimensional centrally symmetric faces of permutation polytopes?
}\end{question}

It should be true that bipyramids over centrally symmetric faces are again realizable as faces of permutation polytopes. 
Even more, we expect that it may be possible to generalize the construction of Theorem \ref{constr}.

\begin{conjecture}{\rm 
The free sum of centrally symmetric faces of permutation polytopes can be combinatorially 
realized as a face of a permutation polytope.
}\end{conjecture}

The next conjecture is based upon explicit checks in low dimensions.

\begin{conjecture}{\rm 
Let $F$ be a face of a permutation polytope. Then the wedge over a face $F'$ of $F$ can be combinatorially 
realized as a face of a permutation polytope, if (or only if) the complement of $F'$ in $F$ is a face of~$F$.
}\end{conjecture}

\subsection{Faces of permutation polytopes given by subgroups}

It would be interesting to know which subgroups of a
permutation group yield faces. One obvious class of such subgroups are
stabilizers. 

For this let us partition $[n] := \{1, \ldots, n\} = \bigsqcup I_i$. Then the polytope of the
stabilizer of this partition $$\operatorname{stab}(G;(I_i)_i) := \{
\sigma \in G : \sigma(I_i)=I_i \text{ for all } i \} \leq G$$ is a face of
$P(G)$. 

\begin{conjecture}{\rm 
    Let $G \le S_n$. Suppose $H \le G$ is a subgroup such that $P(H)
    \preceq P(G)$ is a face. Then $H = \operatorname{stab}(G;(I_i)_i)$
    for a partition $[n] = \bigsqcup I_i$.
    \label{mainconj}}
\end{conjecture}

\begin{prop}
  Conjecture \ref{mainconj} holds for $G = S_n$ and for $G \leq S_n$
  cyclic.
\end{prop}

\begin{proof}
First, let $H \le G = S_n$ with $P(H)$ a face of $P(G) = B_n$. Let
$[n] = \bigsqcup I_i$ be the orbit partition of $H$.
Then $H \le \operatorname{stab}(G;(I_i)_i) = \prod S_{I_i}$. We show
that equality holds.
The face $P(H) \preceq B_n$, is the intersection of the facets
  containing it. For $J = \{ (i,j) : \sigma(i) \neq j \text{ for all }
  \sigma \in H \}$ that means, cf. Remark \ref{strangefaces}, 
  $$P(H) = \conv \{ A \in B_n : a_{ij} = 0 \text{ for all } (i,j) \in
  J \} .$$
  Because $H$ is transitive on $I_i$, we get $J \cap (I_i \times I_i) =
  \emptyset$ so that $H \ge \prod S_{I_i}$.

Second, let $G$ be a cyclic subgroup of $S_n$, and let $H$ be a
subgroup of $G$. Let $G = \< g \>$ and let $g = g_1 \circ \cdots \circ
g_r$ be the cycle decomposition of $g$. Then $H = \< h \>$ for some $h
= g^k = g_1^k \circ \cdots \circ g_r^k$. If the length of $g_i$ is
$z_i$, then $g^k_i$ splits into cycles of length $s_i:=
z_i/\gcd(z_i,k)$. Let $O_1, \ldots , O_l$ be the respective orbits of
$h$. Then $H$ stabilizes the partition  
\[[n] = O_1 \sqcup \cdots
\sqcup O_l \sqcup ([n]\setminus{(O_1 \cup ... \cup O_l)}).\]
Let $h' = g^t$ be another element in $G$ which also stabilizes this
partition. This implies  $z_i/\gcd(z_i,t) \leq s_i$, thus $\gcd(z_i,k)
\leq \gcd(z_i,t)$ for $i = 1, \ldots, r$.
Since $k$ and $t$ divide $\card{G} = o(g) = \lcm(z_1, \ldots, z_r)$,
this yields that $k$ divides $t$. Hence $h' = g^t \in H = \<g^k\>$,
and $H$ is the full stabilizer in $G$ of a partition of $[n]$.
\end{proof}

\begin{acknowledgement}{\rm 
Many of these results are based on extensive calculations using 
the software packages \texttt{GAP} \cite{GAP06} and \texttt{polymake} \cite{GJ05}. We thank Oswin Aichholzer 
for providing the list of $0/1$-polytopes \cite{Aic00}. 
The last three authors were supported by Emmy Noether fellowship HA 4383/1 of
the German Research Foundation (DFG).}
\end{acknowledgement}

\end{document}